\documentclass[12pt,reqno,centertags]{amsart}
\usepackage{amssymb}
\usepackage{amsfonts}
\usepackage{geometry}
\usepackage{fourier}
\newtheorem{theorem}{Theorem}
\theoremstyle{plain}

\newtheorem{proposition}{Proposition}

\numberwithin{equation}{section}

\begin{document}
\title[]{Ruled and quadric surfaces in the 3-dimensional\\ Euclidean space satisfying $\Delta ^{III}\boldsymbol{x}%
=\varLambda \boldsymbol{x}$}
\author{Hassan Al-Zoubi}
\address{Department of Basic Sciences, Al-Zaytoonah University of Jordan,
Amman}
\email{dr.hassanz@zuj.edu.jo}
\author{Stylianos Stamatakis}
\address{Department of Mathematics, Aristotle University of Thessaloniki}
\email{stamata@math.auth.gr}
\date{}

\begin{abstract}
We consider ruled and quadric surfaces in the 3-dimensional Euclidean space
which are of coordinate finite type with respect to the third fundamental
form $III$, i.e., their position vector $\boldsymbol{x}$ satisfies the
relation $\Delta ^{III}\boldsymbol{x}=\varLambda \boldsymbol{x}$ where $\varLambda
$ is a square matrix of order 3. We show that helicoids and spheres are the
only surfaces in $E^3$ satisfying the preceding relation.
\end{abstract}

\maketitle

\section{Introduction}

Let $S$\ be a (connected) surface in a Euclidean 3-space $E^{3}$ referred to
any system of coordinates $u^{1},\ u^{2}$, which does not contain parabolic
points, we denote by $b_{ij}$\ the components of the second fundamental form
$II=b_{ij}du^{i}du^{j}$ of $S.$\ Let $\varphi (u^{1},u^{2})$\ and $\psi
(u^{1},u^{2})$\ be two sufficient differentiable functions on $S$. Then the
first differential parameter of Beltrami with respect to the second
fundamental form\ of $S$\ is defined by

\begin{equation*}
\nabla ^{II}\left( \varphi ,\psi \right) :=b^{ij}\varphi _{/i}\psi _{/j},
\end{equation*}%
\ where $\varphi _{/i}:=\frac{\partial \varphi }{\partial u^{i}}$\ and $%
\left( b^{ij}\right) $\ denotes the inverse tensor of $\left( b_{ij}\right) $%
.

Let $e_{ij}$\ be the components of the third fundamental form $III$ of $S.$
Then the second differential parameter of Beltrami with respect to the third
fundamental form\ of $S$\ is defined by \cite{R9}

\begin{equation*}
\Delta ^{III}\varphi :=-\frac{1}{\sqrt{e}}(\sqrt{e}e^{ij}\varphi _{/i})_{/j}%
\footnote{%
with sign convention such that $\Delta =-\frac{\partial ^{2}}{\partial x^{2}}%
-\frac{\partial ^{2}}{\partial y^{2}}$ for the metric $ds^{2}=dx^{2}$ $%
+dy^{2}.$},
\end{equation*}%
where $\left( e^{ij}\right) $\ denotes the inverse tensor of $\left(
e_{ij}\right) $ and $e:=\det (e_{ij}).$

In \cite{R9}, S. Stamatakis and H. Al-Zoubi showed, for the position vector $%
\boldsymbol{x}=\boldsymbol{x}(u^{1},u^{2})$ of $S$, the relation

\begin{equation}
\Delta ^{III}\boldsymbol{x}=\nabla ^{III}(\frac{2H}{K},\boldsymbol{n})-\frac{%
2H}{K}\boldsymbol{n},  \label{eq 2}
\end{equation}%
where $\boldsymbol{n}$ is the Gauss map, $K$ the Gauss curvature and $H$ the
mean curvature of $S$. Moreover, in this context, the same authors proved
that the surfaces $S:\boldsymbol{x}=\boldsymbol{x}(u^{1},u^{2})$ satisfying
the condition

\begin{equation*}
\Delta ^{III}\boldsymbol{x}=\lambda \boldsymbol{x},\ \ \ \ \lambda \in
\mathbb{R}
,
\end{equation*}%
i.e., for which all coordinate functions are eigenfunctions of $\Delta
^{III} $ with the same eigenvalue $\lambda ,$ are precisely either the
minimal surfaces $\left( \lambda =0\right) ,$ or the part of spheres $\left(
\lambda =2\right) .$

In \cite{R2} B.-Y. Chen introduced the notion of Euclidean immersions of
finite type. In terms of B.-Y. Chen theory, a surface $S$ is said to be of
finite type, if its coordinate functions are a finite sum of eigenfunctions
of\ the Beltrami operator $\Delta ^{III}.$ Therefore the two facts mentioned
above can be stated as follows

\begin{itemize}
\item $S$ is minimal if and only if $S$ is of null type1.

\item $S$ lies in an ordinary sphere S$^{2}$ if and only if $S$ is of type1.
\end{itemize}

Following \cite{R2}\ we say that a surface $S$ is of finite type with
respect to the fundamental form $III$, or briefly of finite $III$-type, if
the position vector $\boldsymbol{x}$ of $S$ can be written as a finite sum\
of nonconstant eigenvectors of the operator $\Delta ^{III},$ i.e., if

\begin{equation}
\boldsymbol{x}=\boldsymbol{x}_{0}+\sum_{i=1}^{m}\boldsymbol{x}_{i},\ \ \ \ \
\ \Delta ^{III}\boldsymbol{x}_{i}=\lambda _{i}\boldsymbol{x}_{i},\ \ \
i=1,...,m,  \label{eq 3}
\end{equation}%
where $\boldsymbol{x}_{0}$ is a fixed vector and $\lambda _{1},\ \lambda
_{2},...,\lambda _{m}$ are eigenvalues of $\Delta ^{III};$ when there are
exactly $k$ nonconstant eigenvectors $\boldsymbol{x}_{1},\ \boldsymbol{x}%
_{2},...,\boldsymbol{x}_{k}$ appearing in (\ref{eq 3}) which all belong to
different eigenvalues $\lambda _{1},\ \lambda _{2},...,\lambda _{k},$ then $%
S $ is said to be of $III$-type $k$, otherwise $S$ is said to be of infinite
type.\textit{\ }When $\lambda _{i}=0$ for some $i=1,2,...,k,$ then $S$ is
said to be of null $III$-type $k$.

Up to now, very little is known about surfaces of\textit{\ }finite $III$%
\textit{-}type\textit{. }Concerning this problem, the only known surfaces of%
\textit{\ }finite $III$\textit{-}type in\textit{\ }$E^{3\text{ }}$are parts
of spheres, the minimal surfaces and the parallel of the minimal surfaces
which are of\textit{\ }null\textit{\ }$III$\textit{-}type 2, (see \cite{R9}).

In this paper we shall be concerned with the ruled and quadric surfaces in $%
E^{3}$ which are connected, complete and which are of coordinate finite $III$%
\textit{-}type, i.e., their position vector $\boldsymbol{x}$ $=\boldsymbol{x}%
(u^{1},u^{2})$ satisfies the relation

\begin{equation}
\Delta ^{III}\boldsymbol{x}=\varLambda \boldsymbol{x},  \label{eq 4}
\end{equation}%
where $\varLambda $ is a square matrix of order 3.

In \cite{R6} F. Dillen, J. Pas and L. Verstraelen studied coordinate finite
type with respect to the first fundamental form $I=g_{ij}du^{i}du^{j}$ and
they proved that

\begin{theorem}
\label{TA} The only surfaces in $%
\mathbb{R}
^{3}$ satisfying%
\begin{equation*}
\Delta ^{I}\boldsymbol{x}=A\boldsymbol{x}+B,\ \ \ \ \ \ A\in M(3,3),\ \ \ \
\ B\in M(3,1)
\end{equation*}%
are the minimal surfaces, the spheres and the circular cylinders.
\end{theorem}

(\textit{M}(m,n) denotes the set of all matrices of the type (m,n)).

While in \cite{R7} O. Garay showed that

\begin{theorem}
\label{TB} The only complete surfaces of revolutions in $%
\mathbb{R}
^{3},$ whose component functions are eigenfunctions of their Laplacian, then
the surface must be a catenoid, a sphere or a right circular cylinder.
\end{theorem}

Recently, H. Al-Zoubi and S. Stamatakis studied coordinate finite type with
respect to the third fundamental form, more precisely, in \cite{R10} they
proved the following

\begin{theorem}
\label{TC} A surface of revolution $S$ in $%
\mathbb{R}
^{3},$ satisfies (\ref{eq 4}), if and only if S is a catenoid or a part of a
sphere.
\end{theorem}

\section{Main results}

Our main results are the following

\begin{proposition}
\label{T2.1} The only \textit{\ ruled surfaces in the 3-dimensional
Euclidean space that satisfies (\ref{eq 4}), are the helicoids.}

\begin{proposition}
\label{T2.2} The only quadric surfaces in the 3-dimensional \textit{%
Euclidean }space that satisfies (\ref{eq 4}), are the spheres.
\end{proposition}
\end{proposition}

Our discussion is local, which means that we show in fact that any open part
of a ruled or a quadric satisfies (\ref{eq 4}), if it is an open part of a
helicoid or\ an open part of a sphere respectively.

Before starting the proof of our main results we first show that the
surfaces mentioned in the above propositions indeed satisfy the condition (%
\ref{eq 4}). On a helicoid the mean curvature vanishes, so, by virtue of (%
\ref{eq 2}), $\Delta ^{III}\boldsymbol{x}=\boldsymbol{0}.$ Therefore a
helicoid satisfies (\ref{eq 4}), where $\varLambda $ is the null matrix in $%
M(3,3)$.

Let S$^{2}(r)$ be a sphere of radius $r$ centered at the origin. If $%
\boldsymbol{x}$ denotes the position vector field of S$^{2}(r),$ then the
Gauss map $\boldsymbol{n}$\ is given by -$\frac{\boldsymbol{x}}{r}.$ For the
Gauss curvature $K$ and the mean curvature $H$ of S$^{2}(r)$ we have $K$ = $%
\frac{1}{r^{2}}$ and $H$ = $\frac{1}{r},$ so, by virtue of (\ref{eq 2}), we
obtain

\begin{equation*}
\Delta ^{III}\boldsymbol{x}=2\boldsymbol{x}
\end{equation*}%
and we find that S$^{2}(r)$ satisfies (\ref{eq 4}) with

\begin{equation*}
\varLambda =\left[
\begin{array}{ccc}
2 & 0 & 0 \\
0 & 2 & 0 \\
0 & 0 & 2%
\end{array}%
\right] .
\end{equation*}

\section{Proof of Proposition \protect\ref{T2.1}}

Let $S$ be a ruled surface in $E^{3}.$ We suppose that $S$ is a
non-cylindrical ruled surface. This surface can be expressed in terms of a
directrix curve $\boldsymbol{\alpha }(s)$ and a unit vectorfield $%
\boldsymbol{\beta }(s)$ pointing along the rulings as

\begin{equation*}
S:\boldsymbol{x}(s,t)=\boldsymbol{\alpha }(s)+t\boldsymbol{\beta }(s),\ \ \
\ \ \ \ \ s\in J,\ \ \ \ \ \ \ \ -\infty <t<\infty .
\end{equation*}%
Moreover, we can take the parameter $s$ to be the arc length along the
spherical curve $\boldsymbol{\beta }(s)$. Then we have

\begin{equation}
\langle \boldsymbol{\alpha }^{\prime },\boldsymbol{\beta }\rangle =0,\ \ \ \
\ \langle \boldsymbol{\beta },\boldsymbol{\beta }\rangle =1,\ \ \ \ \
\langle \boldsymbol{\beta }^{\prime },\boldsymbol{\beta }^{\prime }\rangle
=1,  \label{eq 5}
\end{equation}%
where the prime denotes the derivative in s. The first fundamental form of $%
S $ is

\begin{equation*}
I=qds^{2}+dt^{2},
\end{equation*}%
while the second fundamental form is

\begin{equation*}
II=\frac{p}{\sqrt{q}}ds^{2}+\frac{2A}{\sqrt{q}}dsdt,
\end{equation*}%
where

\begin{eqnarray*}
q &:&=\langle \boldsymbol{\alpha }^{\prime },\boldsymbol{\alpha }^{\prime
}\rangle +2\langle \boldsymbol{\alpha }^{\prime },\boldsymbol{\beta }%
^{\prime }\rangle t+t^{2}, \\
p &:&=\left( \boldsymbol{\alpha }^{\prime },\boldsymbol{\beta },\boldsymbol{%
\alpha }^{\prime \prime }\right) +\left[ \left( \boldsymbol{\alpha }^{\prime
},\boldsymbol{\beta },\boldsymbol{\beta }^{\prime \prime }\right) +\left(
\boldsymbol{\beta }^{\prime },\boldsymbol{\beta },\boldsymbol{\alpha }%
^{\prime \prime }\right) \right] t+\left( \boldsymbol{\beta }^{\prime },%
\boldsymbol{\beta },\boldsymbol{\beta }^{\prime \prime }\right) t^{2}, \\
A &:&=\left( \boldsymbol{\alpha }^{\prime },\boldsymbol{\beta },\boldsymbol{%
\beta }^{\prime }\right) .
\end{eqnarray*}

For convenience, we put

\begin{eqnarray*}
\kappa &:&=\langle \boldsymbol{\alpha }^{\prime },\boldsymbol{\alpha }%
^{\prime }\rangle ,\ \ \ \ \ \ \ \ \ \ \lambda :=\langle \boldsymbol{\alpha }%
^{\prime },\boldsymbol{\beta }^{\prime }\rangle , \\
\mu &:&=\left( \boldsymbol{\beta }^{\prime },\boldsymbol{\beta },\boldsymbol{%
\beta }^{\prime \prime }\right) ,\ \ \ \ \ \nu :=\left( \boldsymbol{\alpha }%
^{\prime },\boldsymbol{\beta },\boldsymbol{\beta }^{\prime \prime }\right)
+\left( \boldsymbol{\beta }^{\prime },\boldsymbol{\beta },\boldsymbol{\alpha
}^{\prime \prime }\right) , \\
\rho &:&=\left( \boldsymbol{\alpha }^{\prime },\boldsymbol{\beta },%
\boldsymbol{\alpha }^{\prime \prime }\right) ,
\end{eqnarray*}%
and thus we have

\begin{equation*}
q=t^{2}+2\lambda t+\kappa ,\ \ \ \ \ \ \ \ \ \ p=\mu t^{2}+\nu t+\rho .
\end{equation*}

For the Gauss curvature $K$ of $S$ we find

\begin{equation}
K=-\frac{A^{2}}{q^{2}}.  \label{eq 6}
\end{equation}

The Beltrami operator with respect to the third fundamental form, after a
lengthy computation, can be expressed as follows

\begin{eqnarray}
\Delta ^{III} &=&-\frac{q}{A^{2}}\frac{\partial ^{2}}{\partial s^{2}}+\frac{%
2qp}{A^{3}}\frac{\partial ^{2}}{\partial s\partial t}-\left( \frac{q^{2}}{%
A^{2}}+\frac{qp^{2}}{A^{4}}\right) \frac{\partial ^{2}}{\partial t^{2}}
\notag \\
&&+\left( \frac{q_{s}}{2A^{2}}+\frac{qp_{t}}{A^{3}}-\frac{pq_{t}}{2A^{3}}%
\right) \frac{\partial }{\partial s}  \notag \\
&&+\left( \frac{qp_{s}}{A^{3}}-\frac{pq_{s}}{2A^{3}}-\frac{pqA^{\prime }}{%
A^{4}}-\frac{qq_{t}}{2A^{2}}+\frac{p^{2}q_{t}}{2A^{4}}-\frac{2qpp_{t}}{A^{4}}%
\right) \frac{\partial }{\partial t}  \notag \\
&=&Q_{1}\frac{\partial ^{2}}{\partial s^{2}}+Q_{2}\frac{\partial ^{2}}{%
\partial s\partial t}+Q_{3}\frac{\partial }{\partial s}+Q_{4}\frac{\partial
}{\partial t}+Q_{5}\frac{\partial ^{2}}{\partial t^{2}}  \label{eq 7}
\end{eqnarray}%
where

\begin{equation*}
q_{t}:=\frac{\partial q}{\partial t},\ \ \ \ q_{s}:=\frac{\partial q}{%
\partial s},\ \ \ \ p_{t}:=\frac{\partial p}{\partial t},\ \ \ \ p_{s}:=%
\frac{\partial p}{\partial s}
\end{equation*}%
and $Q_{1},Q_{2},...,Q_{5}\ $are polynomials in $t$ with functions in $s$ as
coefficients, and $deg(Q_{i})\leq 6.$ More precisely we have

\begin{equation*}
Q_{1}=-\frac{1}{A^{2}}[t^{2}+2\lambda t+\kappa ],
\end{equation*}%
\begin{eqnarray*}
Q_{2} &=&\frac{2}{A^{3}}[\mu t^{4}+\left( 2\lambda \mu +\nu \right)
t^{3}+\left( 2\lambda \nu +\rho +\kappa \mu \right) t^{2} \\
&&+\left( 2\lambda \rho +\kappa \nu \right) t+\kappa \rho ],
\end{eqnarray*}%
\begin{eqnarray*}
Q_{3} &=&\frac{1}{A^{3}}[\mu t^{3}+3\lambda \mu t^{2}+\left( \lambda \nu
-\rho +2\kappa \mu +\lambda ^{\prime }A\right) t \\
&&+\frac{1}{2}\kappa ^{\prime }A-\lambda \rho +\kappa \nu ],
\end{eqnarray*}%
\begin{eqnarray*}
Q_{4} &=&\frac{1}{A^{4}}[-3\mu ^{2}t^{5}+\left( \mu ^{\prime }A-\mu
A^{\prime }-4\mu \nu -7\lambda \mu ^{2}\right) t^{4} \\
&&+(\nu ^{\prime }A-\nu A^{\prime }+2\lambda \mu ^{\prime }A-2\lambda \mu
A^{\prime }-\lambda ^{\prime }\mu A-A^{2} \\
&&-10\lambda \mu \nu -2\mu \rho -\nu ^{2}-4\kappa \mu ^{2})t^{3} \\
&&+(\kappa \mu ^{\prime }A-\kappa \mu A^{\prime }-\frac{1}{2}\kappa ^{\prime
}\mu A+2\lambda \nu ^{\prime }A-2\lambda \nu A^{\prime }-\lambda ^{\prime
}\nu A \\
&&-\rho A^{\prime }+\rho ^{\prime }A-3\lambda A^{2}-3\lambda \nu
^{2}-6\lambda \mu \rho -6\kappa \mu \nu )t^{2} \\
&&+(\kappa \nu ^{\prime }A-\kappa \nu A^{\prime }-\frac{1}{2}\kappa ^{\prime
}\nu A+2\lambda \rho ^{\prime }A-2\lambda \rho A^{\prime }-\lambda ^{\prime
}\rho A \\
&&-\kappa A^{2}-2\lambda ^{2}A^{2}-2\kappa \nu ^{2}+\rho ^{2}-2\lambda \nu
\rho -4\kappa \mu \rho )t \\
&&+(\kappa \rho ^{\prime }A-\kappa \rho A^{\prime }-\frac{1}{2}\kappa
^{\prime }\rho A+\lambda \rho ^{2}-\kappa \lambda A^{2}-2\kappa \nu \rho )],
\end{eqnarray*}%
\begin{eqnarray*}
Q_{5} &=&-\frac{1}{A^{4}}[\mu ^{2}t^{6}+\left( 2\mu \nu +2\lambda \mu
^{2}\right) t^{5} \\
&&+(2\mu \rho +\nu ^{2}+4\lambda \mu \nu +\kappa \mu ^{2}+A^{2})t^{4} \\
&&+(2\nu \rho +4\lambda \mu \rho +2\lambda \nu ^{2}+2\kappa \mu \nu
+4\lambda A^{2})t^{3} \\
&&+(\rho ^{2}+4\lambda \nu \rho +2\kappa \mu \rho +\kappa \nu ^{2}+4\lambda
^{2}A^{2}+2\kappa A^{2})t^{2} \\
&&+(2\lambda \rho ^{2}+2\kappa \nu \rho +4\lambda \kappa A^{2})t+(\kappa
\rho ^{2}+\kappa ^{2}A^{2})].
\end{eqnarray*}

Applying (\ref{eq 7}) for the position vector $\boldsymbol{x}$ one finds%
\begin{equation}
\Delta ^{III}\boldsymbol{x}=Q_{1}\boldsymbol{\alpha }^{\prime \prime }+Q_{2}%
\boldsymbol{\beta }^{\prime }+Q_{3}\boldsymbol{\alpha }^{\prime }+Q_{4}%
\boldsymbol{\beta }+\left( Q_{1}\boldsymbol{\beta }^{\prime \prime }+Q_{3}%
\boldsymbol{\beta }^{\prime }\right) t.  \label{eq 8}
\end{equation}

Let $\boldsymbol{x}=(x_{1},x_{2},x_{3}),\boldsymbol{\alpha }=(\alpha
_{1},\alpha _{2},\alpha _{3})$ and $\boldsymbol{\beta }=(\beta _{1},\beta
_{2},\beta _{3})$ be the coordinate functions of $\boldsymbol{x},\boldsymbol{%
\alpha }$ and $\boldsymbol{\beta }.$ By virtue of (\ref{eq 8}) we obtain%
\begin{eqnarray}
\Delta ^{III}x_{i} &=&Q_{1}\alpha _{i}^{\prime \prime }+Q_{2}\beta
_{i}^{\prime }+Q_{3}\alpha _{i}^{\prime }+Q_{4}\beta _{i}+\left( Q_{1}\beta
_{i}^{\prime \prime }+Q_{3}\beta _{i}^{\prime }\right) t,  \label{eq 9} \\
i &=&1,2,3.  \notag
\end{eqnarray}

We denote by $\lambda _{ij},i,j=1,2,3$ the entries of the matrix $\varLambda $.
Using (\ref{eq 9}) and condition (\ref{eq 4}) we have

\begin{eqnarray*}
&&Q_{1}\alpha _{i}^{\prime \prime }+Q_{2}\beta _{i}^{\prime }+Q_{3}\alpha
_{i}^{\prime }+Q_{4}\beta _{i}+\left( Q_{1}\beta _{i}^{\prime \prime
}+Q_{3}\beta _{i}^{\prime }\right) t \\
&=&\lambda _{i1}\alpha _{1}+\lambda _{i2}\alpha _{2}+\lambda _{i3}\alpha
_{3}+(\lambda _{i1}\beta _{1}+\lambda _{i2}\beta _{2}+\lambda _{i3}\beta
_{3})t, \\
i &=&1,2,3.
\end{eqnarray*}

Consequently%
\begin{eqnarray}
&&-3\mu ^{2}\beta _{i}t^{5}+[\left( \mu ^{\prime }A-\mu A^{\prime }-4\mu \nu
-7\lambda \mu ^{2}\right) \beta _{i}+3\mu A\beta _{i}^{\prime }]t^{4}  \notag
\\
&&+[\mu A\alpha _{i}^{\prime }-A^{2}\beta _{i}^{\prime \prime }+(2\nu
A+7\lambda \mu A)\beta _{i}^{\prime }+(\nu ^{\prime }A-\nu A^{\prime }
\notag \\
&&+2\lambda \mu ^{\prime }A-2\lambda \mu A^{\prime }-\lambda ^{\prime }\mu
A-A^{2}-10\lambda \mu \nu -2\mu \rho -\nu ^{2}  \notag \\
&&-4\kappa \mu ^{2})\beta _{i}]t^{3}+[(\kappa \mu ^{\prime }A-\kappa \mu
A^{\prime }-\frac{1}{2}\kappa ^{\prime }\mu A+2\lambda \nu ^{\prime }A
\notag \\
&&-2\lambda \nu A^{\prime }-\lambda ^{\prime }\nu A-\rho A^{\prime }+\rho
^{\prime }A-3\lambda A^{2}-3\lambda \nu ^{2}-6\lambda \mu \rho   \notag \\
&&-6\kappa \mu \nu )\beta _{i}+3\lambda \mu A\alpha _{i}^{\prime }-2\lambda
A^{2}\beta _{i}^{\prime \prime }-A^{2}\alpha _{i}^{\prime \prime }+(\lambda
^{\prime }A+5\lambda \nu   \notag \\
&&+4\kappa \mu +\rho )A\beta _{i}^{\prime }]t^{2}+[(\frac{1}{2}\kappa
^{\prime }A+3\kappa \nu +3\lambda \rho )A\beta _{i}^{\prime }  \notag \\
&&+(\kappa \nu ^{\prime }A-\kappa \nu A^{\prime }-\frac{1}{2}\kappa ^{\prime
}\nu A+2\lambda \rho ^{\prime }A-2\lambda \rho A^{\prime }-\lambda ^{\prime
}\rho A  \notag \\
&&-\kappa A^{2}-2\lambda ^{2}A^{2}-2\kappa \nu ^{2}+\rho ^{2}-2\lambda \nu
\rho -4\kappa \mu \rho )\beta _{i}  \notag \\
&&-2\lambda A^{2}\alpha _{i}^{\prime \prime }-\kappa A^{2}\beta _{i}^{\prime
\prime }+\left( \lambda \nu -\rho +2\kappa \mu +\lambda ^{\prime }A\right)
A\alpha _{i}^{\prime }]t  \notag \\
&&-A^{4}\left( \lambda _{i1}\beta _{1}+\lambda _{i2}\beta _{2}+\lambda
_{i3}\beta _{3}\right) t+(\kappa \rho ^{\prime }A-\kappa \rho A^{\prime }
\notag \\
&&-\frac{1}{2}\kappa ^{\prime }\rho A+\lambda \rho ^{2}-\kappa \lambda
A^{2}-2\kappa \nu \rho )\beta _{i}-\kappa A^{2}\alpha _{i}^{\prime \prime }
\notag \\
&&+2\kappa \rho A\beta _{i}^{\prime }+(\frac{1}{2}\kappa ^{\prime }A-\lambda
\rho +\kappa \nu )A\alpha _{i}^{\prime }  \notag \\
&&-A^{4}\left( \lambda _{i1}\alpha _{1}+\lambda _{i2}\alpha _{2}+\lambda
_{i3}\alpha _{3}\right) =0.  \label{eq 10}
\end{eqnarray}

For $i=1,2,3,$ (\ref{eq 10}) is a polynomial in $t$ with functions in $s$ as
coefficients. This implies that the coefficients of the powers of $t$ in (%
\ref{eq 10}) must be zeros, so we obtain, for $i=1,2,3,$ the following
equations

\begin{equation}
3\mu ^{2}\beta _{i}=0,  \label{eq 11}
\end{equation}

\begin{equation*}
\left( \mu ^{\prime }A-\mu A^{\prime }-4\mu \nu -7\lambda \mu ^{2}\right)
\beta _{i}+3\mu A\beta _{i}^{\prime }=0,
\end{equation*}%
\begin{eqnarray}
&&\mu A\alpha _{i}^{\prime }-A^{2}\beta _{i}^{\prime \prime }+(2\nu
A+7\lambda \mu A)\beta _{i}^{\prime }  \notag \\
&&+(\nu ^{\prime }A-\nu A^{\prime }+2\lambda \mu ^{\prime }A-2\lambda \mu
A^{\prime }-\lambda ^{\prime }\mu A  \notag \\
&&-A^{2}-10\lambda \mu \nu -2\mu \rho -\nu ^{2}-4\kappa \mu ^{2})\beta
_{i}=0,  \label{eq 12}
\end{eqnarray}%
\begin{eqnarray}
&&(\kappa \mu ^{\prime }A-\kappa \mu A^{\prime }-\frac{1}{2}\kappa ^{\prime
}\mu A+2\lambda \nu ^{\prime }A-2\lambda \nu A^{\prime }-\lambda ^{\prime
}\nu A  \notag \\
&&-\rho A^{\prime }+\rho ^{\prime }A-3\lambda A^{2}-3\lambda \nu
^{2}-6\lambda \mu \rho -6\kappa \mu \nu )\beta _{i}+3\lambda \mu A\alpha
_{i}^{\prime }  \notag \\
&&-2\lambda A^{2}\beta _{i}^{\prime \prime }-A^{2}\alpha _{i}^{\prime \prime
}+(\lambda ^{\prime }A+5\lambda \nu +4\kappa \mu +\rho )A\beta _{i}^{\prime
}=0,  \label{eq 13}
\end{eqnarray}%
\begin{eqnarray}
&&(\kappa \nu ^{\prime }A-\kappa \nu A^{\prime }-\frac{1}{2}\kappa ^{\prime
}\nu A+2\lambda \rho ^{\prime }A-2\lambda \rho A^{\prime }-\lambda ^{\prime
}\rho A-\kappa A^{2}  \notag \\
&&-2\lambda ^{2}A^{2}-2\kappa \nu ^{2}+\rho ^{2}-2\lambda \nu \rho -4\kappa
\mu \rho )\beta _{i}-2\lambda A^{2}\alpha _{i}^{\prime \prime }-\kappa
A^{2}\beta _{i}^{\prime \prime }  \notag \\
&&+(\frac{1}{2}\kappa ^{\prime }A+2\kappa \nu +4\lambda \rho )A\beta
_{i}^{\prime }+\left( \lambda \nu -\rho +2\kappa \mu +\lambda ^{\prime
}A\right) A\alpha _{i}^{\prime }  \notag \\
&=&A^{4}(\lambda _{i1}\beta _{1}+\lambda _{i2}\beta _{2}+\lambda _{i3}\beta
_{3}),  \label{eq 14}
\end{eqnarray}%
\begin{eqnarray}
&&(\kappa \rho ^{\prime }A-\kappa \rho A^{\prime }-\frac{1}{2}\kappa
^{\prime }\rho A+\lambda \rho ^{2}-\kappa \lambda A^{2}-2\kappa \nu \rho
)\beta _{i}+2\kappa \rho A\beta _{i}^{\prime }  \notag \\
&&+(\frac{1}{2}\kappa ^{\prime }A-\lambda \rho +\kappa \nu )A\alpha
_{i}^{\prime }-\kappa A^{2}\alpha _{i}^{\prime \prime }  \notag \\
&=&A^{4}\left( \lambda _{i1}\alpha _{1}+\lambda _{i2}\alpha _{2}+\lambda
_{i3}\alpha _{3}\right) .  \label{eq 15}
\end{eqnarray}

From (\ref{eq 11}) one finds

\begin{equation}
\mu =\left( \boldsymbol{\beta }^{\prime },\boldsymbol{\beta },\boldsymbol{%
\beta }^{\prime \prime }\right) =0,  \label{eq 16}
\end{equation}%
which implies that the vectors $\boldsymbol{\beta }^{\prime },\boldsymbol{%
\beta },\boldsymbol{\beta }^{\prime \prime }$ are linearly dependent, and
hence there exist two functions $\sigma _{1}=\sigma _{1}(s)$ and $\sigma
_{2}=\sigma _{2}(s)$ such that

\begin{equation}
\boldsymbol{\beta }^{\prime \prime }=\sigma _{1}\boldsymbol{\beta +}\sigma
_{2}\boldsymbol{\beta }^{\prime }.  \label{eq 17}
\end{equation}

On differentiating $\langle \boldsymbol{\beta }^{\prime },\boldsymbol{\beta }%
^{\prime }\rangle =1,$ we obtain $\langle \boldsymbol{\beta }^{\prime },%
\boldsymbol{\beta }^{\prime \prime }\rangle =0.$ So from (\ref{eq 17}) we
have

\begin{equation}
\boldsymbol{\beta }^{\prime \prime }=\sigma _{1}\boldsymbol{\beta }.
\label{eq 19}
\end{equation}

By taking the derivative of $\langle \boldsymbol{\beta },\boldsymbol{\beta }%
\rangle =1$ twice, we find that

\begin{equation*}
\langle \boldsymbol{\beta }^{\prime },\boldsymbol{\beta }^{\prime }\rangle
+\langle \boldsymbol{\beta },\boldsymbol{\beta }^{\prime \prime }\rangle =0.
\end{equation*}

But $\langle \boldsymbol{\beta }^{\prime },\boldsymbol{\beta }^{\prime
}\rangle =1,$ and taking into account (\ref{eq 19}) we find that $\sigma
_{1}(s)=-1.$ Thus (\ref{eq 19}) becomes $\boldsymbol{\beta }^{\prime \prime
}=-\boldsymbol{\beta }$ which implies that

\begin{equation}
\beta _{i}^{\prime \prime }=-\beta _{i},\ \ \ \ i=1,2,3.  \label{eq 20}
\end{equation}

Using (\ref{eq 16})\ and (\ref{eq 20})\ equation (\ref{eq 12})\ reduces to

\begin{equation*}
2\nu A\beta _{i}^{\prime }+(\nu ^{\prime }A-\nu A^{\prime }-\nu ^{2})\beta
_{i}=0,\ \ \ \ i=1,2,3
\end{equation*}%
or, in vector notation

\begin{equation}
2\nu A\boldsymbol{\beta }^{\prime }+(\nu ^{\prime }A-\nu A^{\prime }-\nu
^{2})\boldsymbol{\beta }=\boldsymbol{0}.  \label{eq 21}
\end{equation}

By taking the derivative of $\langle \boldsymbol{\beta },\boldsymbol{\beta }%
\rangle =1$, we find that the vectors $\boldsymbol{\beta },\boldsymbol{\beta
}^{\prime }$ are linearly independent, and so from (\ref{eq 21}) we obtain
that $\nu A=0.$ We note that $A\neq 0$, since from (\ref{eq 6}) the Gauss
curvature vanishes, so we are left with $\nu =0$. Then equation (\ref{eq 13}%
) becomes

\begin{equation*}
-A^{2}\alpha _{i}^{\prime \prime }+(\lambda ^{\prime }A+\rho )A\beta
_{i}^{\prime }+(\rho ^{\prime }A-\rho A^{\prime }-\lambda A^{2})\beta
_{i}=0,\ \ \ \ i=1,2,3
\end{equation*}%
or, in vector notation

\begin{equation}
-A^{2}\boldsymbol{\alpha }^{\prime \prime }+(\lambda ^{\prime }A+\rho )A%
\boldsymbol{\beta }^{\prime }+(\rho ^{\prime }A-\rho A^{\prime }-\lambda
A^{2})\boldsymbol{\beta }=\boldsymbol{0}.  \label{eq 22}
\end{equation}

Taking the inner product of both sides of the above equation with $%
\boldsymbol{\beta }^{\prime }$ we find in view of (\ref{eq 5}) that

\begin{equation}
-A^{2}\langle \boldsymbol{\alpha }^{\prime \prime },\boldsymbol{\beta }%
^{\prime }\rangle +\rho A+\lambda ^{\prime }A^{2}=0.  \label{eq 23}
\end{equation}

On differentiating $\lambda =\langle \boldsymbol{\alpha }^{\prime },%
\boldsymbol{\beta }^{\prime }\rangle $ with respect to $s,$ by virtue of (%
\ref{eq 20}) and (\ref{eq 5}), we get

\begin{equation}
\lambda ^{\prime }=\langle \boldsymbol{\alpha }^{\prime \prime },\boldsymbol{%
\beta }^{\prime }\rangle +\langle \boldsymbol{\alpha }^{\prime },\boldsymbol{%
\beta }^{\prime \prime }\rangle =\langle \boldsymbol{\alpha }^{\prime \prime
},\boldsymbol{\beta }^{\prime }\rangle -\langle \boldsymbol{\alpha }^{\prime
},\boldsymbol{\beta }\rangle =\langle \boldsymbol{\alpha }^{\prime \prime },%
\boldsymbol{\beta }^{\prime }\rangle .  \label{eq 24}
\end{equation}

Hence, (\ref{eq 23}) reduces to $\rho A=0,$ which implies that $\rho =0.$
Thus the vectors $\boldsymbol{\alpha }^{\prime },\boldsymbol{\beta },%
\boldsymbol{\alpha }^{\prime \prime }$ are linearly dependent, and so there
exist two functions $\sigma _{3}=\sigma _{3}(s)$ and $\sigma _{4}=\sigma
_{4}(s)$ such that

\begin{equation}
\boldsymbol{\alpha }^{\prime \prime }=\sigma _{3}\boldsymbol{\beta +}\sigma
_{4}\boldsymbol{\alpha }^{\prime }.  \label{eq 25}
\end{equation}

Taking the inner product of both sides of the last equation with $%
\boldsymbol{\beta }$ we find in view of (\ref{eq 5}) that $\sigma
_{3}=\langle \boldsymbol{\alpha }^{\prime \prime },\boldsymbol{\beta }%
\rangle .$

Now, by taking the derivative of $\langle \boldsymbol{\alpha }^{\prime },%
\boldsymbol{\beta }\rangle =0,$ we find $\langle \boldsymbol{\alpha }%
^{\prime \prime },\boldsymbol{\beta }\rangle +\langle \boldsymbol{\alpha }%
^{\prime },\boldsymbol{\beta }^{\prime }\rangle =0,$ that is

\begin{equation}
\langle \boldsymbol{\alpha }^{\prime \prime },\boldsymbol{\beta }\rangle
+\lambda =0,  \label{eq 26}
\end{equation}%
and hence $\sigma _{3}=-\lambda .$

Taking again the inner product of both sides of equation (\ref{eq 25}) with $%
\boldsymbol{\beta }^{\prime }$ we find in view of (\ref{eq 5}) that

\begin{equation}
\langle \boldsymbol{\alpha }^{\prime \prime },\boldsymbol{\beta }^{\prime
}\rangle =\sigma _{4}\lambda .  \label{eq 27}
\end{equation}

Using (\ref{eq 24}) we find $\lambda ^{\prime }=\sigma _{4}\lambda .$ Thus $%
\sigma _{4}=\frac{\lambda ^{\prime }}{\lambda }.$ Therefore

\begin{equation}
\boldsymbol{\alpha }^{\prime \prime }=-\lambda \boldsymbol{\beta +}\frac{%
\lambda ^{\prime }}{\lambda }\boldsymbol{\alpha }^{\prime }.  \label{eq 28}
\end{equation}

We distinguish two cases

\textbf{Case 1}: $\lambda =0.$ Because of $\rho =0$ equation (\ref{eq 22})
would yield $A=0,$ which is clearly impossible for the surfaces under
consideration.

\textbf{Case 2}: $\lambda \neq 0.$ From (\ref{eq 22}), (\ref{eq 28}) and $%
\rho =0$ we find that

\begin{equation*}
-\frac{\lambda ^{\prime }}{\lambda }A^{2}\boldsymbol{\alpha }^{\prime
}+\lambda ^{\prime }A^{2}\boldsymbol{\beta }^{\prime }=\mathbf{0}
\end{equation*}%
which implies that $\lambda ^{\prime }(\boldsymbol{\alpha }^{\prime
}-\lambda \boldsymbol{\beta }^{\prime })=\mathbf{0}.$

If $\lambda ^{\prime }\neq 0,$ then $\boldsymbol{\alpha }^{\prime }=\lambda
\boldsymbol{\beta }^{\prime }.$ Hence $\boldsymbol{\alpha }^{\prime },%
\boldsymbol{\beta }^{\prime }$ are linearly dependent, and so $A=0$ which
contradicts our previous assumption. Thus $\lambda ^{\prime }=0.$ From (\ref%
{eq 28}) we have

\begin{equation}
\boldsymbol{\alpha }^{\prime \prime }=-\lambda \boldsymbol{\beta }.
\label{eq 29}
\end{equation}

On the other hand, by taking the derivative of $\kappa $ and using the last
equation we obtain that $\kappa $ is constant. Hence equations (\ref{eq 14})
and (\ref{eq 15}) reduce to

\begin{eqnarray*}
\lambda _{i1}\beta _{1}+\lambda _{i2}\beta _{2}+\lambda _{i3}\beta _{3} &=&0,
\\
\lambda _{i1}\alpha _{1}+\lambda _{i2}\alpha _{2}+\lambda _{i3}\alpha _{3}
&=&0,\ \ i=1,2,3\
\end{eqnarray*}%
and so $\lambda _{ij}=0,i,j=1,2,3.$

Since the parameter $s$ is the arc length of the spherical curve $%
\boldsymbol{\beta }(s)$, and because of (\ref{eq 16}) we suppose, without
loss of generality, that the parametrization of $\boldsymbol{\beta }(s)$ is

\begin{equation*}
\boldsymbol{\beta }(s)=(\cos s,\sin s,0).
\end{equation*}

Integrating (\ref{eq 29}) twice we get

\begin{equation*}
\boldsymbol{\alpha }(s)=(c_{1}s+c_{2}+\lambda \cos s,c_{3}s+c_{4}+\lambda
\sin s,c_{5}s+c_{6}),
\end{equation*}%
where $c_{i},i=1,2,...,6$ are constants.

Since $\kappa =\langle \boldsymbol{\alpha }^{\prime },\boldsymbol{\alpha }%
^{\prime }\rangle $ is constant, it's easy to show that $c_{1}=c_{3}=0.$
Hence $\boldsymbol{\alpha }(s)$ reduces to

\begin{equation*}
\boldsymbol{\alpha }(s)=(c_{2}+\lambda \cos s,c_{4}+\lambda \sin
s,c_{5}s+c_{6}).
\end{equation*}

Thus we have

\begin{equation*}
S:\boldsymbol{x}(s,t)=(c_{2}+(\lambda +t)\cos s,c_{4}+(\lambda +t)\sin
s,c_{5}s+c_{6})
\end{equation*}%
which is a helicoid.

\section{\protect\bigskip Proof of Proposition \protect\ref{T2.2}}

Let now $S$ be a quadric surface in the Euclidean 3-space $E^{3}.$ Then $S$
is either ruled, or of one of the following two kinds

\begin{equation}
z^{2}-ax^{2}-by^{2}=c,\ \ \ \ abc\neq 0  \label{I}
\end{equation}%
or

\begin{equation}
z=\frac{a}{2}x^{2}+\frac{b}{2}y^{2},\ \ \ \ a>0,\ b>0.  \label{II}
\end{equation}

If $S$ is ruled and satisfies (\ref{eq 4}), then by Proposition \ref{T2.1} $%
S $ is a helicoid. We first show that a quadric of the kind (\ref{I})
satisfies (\ref{eq 4}) if and only if $a=-1$ and $b=-1$, which means that $S$
is a sphere. Next we show that a quadric of the kind (\ref{II}) is never
satisfying (\ref{eq 4}).

\subsection{Quadrics of the first kind}

\noindent This kind of quadric surfaces can be parametrized as follows

\begin{equation*}
\boldsymbol{x}(u,v)=\left( u,v,\sqrt{c+au^{2}+bv^{2}}\right) .
\end{equation*}

Let's denote the function $c+au^{2}+bv^{2}$ by $\omega $ and the function $%
c+a(a+1)u^{2}+b(b+1)v^{2}$ by $T.$ Then the components $g_{ij},b_{ij}$ and $%
e_{ij}$ of the first, second and third fundamental tensors in (local)
coordinates are the following

\begin{equation*}
g_{11}=1+\frac{\left( au\right) ^{2}}{\omega },\ \ \ g_{12}=\frac{abuv}{%
\omega },\ \ \ g_{22}=1+\frac{\left( bv\right) ^{2}}{\omega },
\end{equation*}

\begin{equation*}
b_{11}=\frac{a\left( c+bv^{2}\right) }{\omega \sqrt{T}},\ \ \ b_{12}=-\frac{%
abuv}{\omega \sqrt{T}},\ \ \ b_{22}=\frac{b\left( c+au^{2}\right) }{\omega
\sqrt{T}}
\end{equation*}%
and

\begin{eqnarray*}
e_{11} &=&\frac{a^{2}}{\omega T^{2}}\left[
(buv)^{2}+(bv^{2}+c)^{2}+b^{2}v^{2}\omega \right] , \\
e_{12} &=&\frac{ab}{\omega T^{2}}\left[ c(a+b)uv+abuv(u^{2}+v^{2}+\omega )%
\right] , \\
e_{22} &=&\frac{b^{2}}{\omega T^{2}}\left[
(auv)^{2}+(au^{2}+c)^{2}+a^{2}u^{2}\omega \right] .
\end{eqnarray*}

Notice that $\omega $ and $T$ are polynomials in $u$ and $v$. If for
simplicity we put

\begin{equation*}
C(u,v)=(buv)^{2}+(bv^{2}+c)^{2}+b^{2}v^{2}\omega ,
\end{equation*}

\begin{equation*}
B(u,v)=uv\left[ c(a+b)+ab(u^{2}+v^{2}+\omega )\right] ,
\end{equation*}

\begin{equation*}
A(u,v)=(auv)^{2}+(au^{2}+c)^{2}+a^{2}u^{2}\omega ,
\end{equation*}

then the third fundamental tensors $e_{ij}$ turns into

\begin{equation*}
e_{11}=\frac{a^{2}}{\omega T^{2}}C(u,v),\ \ \ e_{12}=-\frac{ab}{\omega T^{2}}%
B(u,v),\ \ \ e_{22}=\frac{b^{2}}{\omega T^{2}}A(u,v).
\end{equation*}

Hence the Beltrami operator $\Delta ^{III}$ of $S$ can be expressed as
follows

\begin{eqnarray}
\Delta ^{III} &=&-\frac{T}{\left( abc\right) ^{2}}\left[ b^{2}A\frac{%
\partial ^{2}}{\partial u^{2}}+2abB\frac{\partial ^{2}}{\partial u\partial v}%
+a^{2}C\frac{\partial ^{2}}{\partial v^{2}}\right]  \notag \\
&&-\frac{T}{\left( abc\right) ^{2}}\left[ b\left( b\frac{\partial A}{%
\partial u}+a\frac{\partial B}{\partial v}\right) \frac{\partial }{\partial u%
}+a\left( a\frac{\partial C}{\partial v}+b\frac{\partial B}{\partial u}%
\right) \frac{\partial }{\partial v}\right]  \notag \\
&&+\frac{T}{\left( abc\right) ^{2}}\left[ \frac{ab^{2}}{\omega }\left(
uA+vB\right) \frac{\partial }{\partial u}+\frac{a^{2}b}{\omega }\left(
uB+vC\right) \frac{\partial }{\partial v}\right]  \notag \\
&&+\frac{1}{\left( abc\right) ^{2}}[ab^{2}\left( \left( a+1\right) uA+\left(
b+1\right) vB\right) \frac{\partial }{\partial u}  \notag \\
&&+a^{2}b\left( \left( b+1\right) vC+\left( a+1\right) uB\right) \frac{%
\partial }{\partial v}].  \label{eq 31}
\end{eqnarray}

We remark that

\begin{equation*}
b\frac{\partial A}{\partial u}+a\frac{\partial B}{\partial v}=au\left[
5ab(a+1)u^{2}+5ab(b+1)v^{2}+c(3ab+5b+a)\right] ,
\end{equation*}

\begin{equation*}
a\frac{\partial C}{\partial v}+b\frac{\partial B}{\partial u}=av\left[
5ab(a+1)u^{2}+5ab(b+1)v^{2}+c(3ab+5a+b)\right] ,
\end{equation*}

\begin{equation*}
uA+vB=\left[ c+a(a+1)u^{2}+a(b+1)v^{2}\right] u\omega ,
\end{equation*}

\begin{equation*}
uB+vC=\left[ c+b(a+1)u^{2}+b(b+1)v^{2}\right] v\omega ,
\end{equation*}

\begin{equation*}
\left( a+1\right) uA+\left( b+1\right) vB=u\left[
c(a+1)+a(a+1)u^{2}+a(b+1)v^{2}\right] T,
\end{equation*}

\begin{equation*}
\left( b+1\right) vC+\left( a+1\right) uB=v\left[
c(b+1)+b(a+1)u^{2}+b(b+1)v^{2}\right] T.
\end{equation*}

We denote by $\lambda _{ij},i,j=1,2,3$ the entries of the matrix $\varLambda $.
On account of (\ref{eq 4}) we get

\begin{equation}
\Delta ^{III}x_{1}=\Delta ^{III}u=\lambda _{11}u+\lambda _{12}v+\lambda _{13}%
\sqrt{\omega },  \label{eq 32}
\end{equation}

\begin{equation}
\Delta ^{III}x_{2}=\Delta ^{III}v=\lambda _{21}u+\lambda _{22}v+\lambda _{23}%
\sqrt{\omega },  \label{eq 33}
\end{equation}

\begin{equation*}
\Delta ^{III}x_{3}=\Delta ^{III}\sqrt{\omega }=\lambda _{31}u+\lambda
_{32}v+\lambda _{33}\sqrt{\omega }.
\end{equation*}

Applying (\ref{eq 31}) on the coordinate functions $x_{i},i=1,2$ of the
position vector $\boldsymbol{x}$ and by virtue of (\ref{eq 32}) and (\ref{eq
33}), we find respectively

\begin{eqnarray}
\Delta ^{III}u &=&-\frac{uT}{c^{2}}\left[ 3\left( a+1\right) u^{2}+3\left(
b+1\right) v^{2}+\frac{c(3b+a+2ab)}{ab}\right]  \notag \\
&=&\lambda _{11}u+\lambda _{12}v+\lambda _{13}\sqrt{\omega },  \label{eq 34}
\end{eqnarray}

\begin{eqnarray}
\Delta ^{III}v &=&-\frac{vT}{c^{2}}\left[ 3\left( a+1\right) u^{2}+3\left(
b+1\right) v^{2}+\frac{c(b+3a+2ab)}{ab}\right]  \notag \\
&=&\lambda _{21}u+\lambda _{22}v+\lambda _{23}\sqrt{\omega },  \label{eq 35}
\end{eqnarray}

Putting $v=0$ in (\ref{eq 34}), we obtain that%
\begin{eqnarray*}
&&-\frac{3a(a+1)^{2}}{c^{2}}u^{5}-\frac{(a+1)(6b+a+2ab)}{bc}u^{3}-\frac{%
\left( 3b+a+2ab\right) }{ab}u \\
&=&\lambda _{11}u+\lambda _{13}\sqrt{c+au^{2}}.
\end{eqnarray*}

Since $a\neq 0$ and $c\neq 0$ this implies that $a=-1.$

Similarly, if we put $u=0$ in (\ref{eq 35}) we obtain that%
\begin{eqnarray*}
&&-\frac{3b(b+1)^{2}}{c^{2}}v^{5}-\frac{(b+1)(b+6a+2ab)}{ac}v^{3}-\frac{%
\left( b+3a+2ab\right) }{ab}v \\
&=&\lambda _{22}v+\lambda _{23}\sqrt{c+bv^{2}}.
\end{eqnarray*}

This implies that $b=-1.$ Hence $S$ must be a sphere.

\subsection{Quadrics of the second kind}

For this kind of surfaces we can consider a parametrization

\begin{equation*}
\boldsymbol{x}(u,v)=\left( u,v,\frac{a}{2}u^{2}+\frac{b}{2}v^{2}\right) .
\end{equation*}

Then the components $g_{ij},b_{ij}$ and $e_{ij}$ of the first, second and
third fundamental tensors are the following

\begin{equation*}
g_{11}=1+\left( au\right) ^{2},\ \ \ g_{12}=abuv,\ \ \ g_{22}=1+\left(
bv\right) ^{2},
\end{equation*}

\begin{equation*}
b_{11}=\frac{a}{\sqrt{g}},\ \ \ b_{12}=0,\ \ \ b_{22}=\frac{b}{\sqrt{g}},
\end{equation*}

\begin{equation*}
e_{11}=\frac{a^{2}}{g^{2}}(1+b^{2}v^{2}),\ \ \ e_{12}=-\frac{a^{2}b^{2}}{%
g^{2}}uv,\ \ \ e_{22}=\frac{b^{2}}{g^{2}}(1+a^{2}u^{2}),
\end{equation*}%
where $g:=\det \left( g_{ij}\right) =1+\left( au\right) ^{2}+\left(
bv\right) ^{2}.$

A straightforward computation shows that the Beltrami operator $\Delta
^{III} $ of $S$ takes the following form

\begin{eqnarray}
\Delta ^{III} &=&-\frac{g(1+a^{2}u^{2})}{a^{2}}\frac{\partial ^{2}}{\partial
u^{2}}-\frac{g(1+b^{2}v^{2})}{b^{2}}\frac{\partial ^{2}}{\partial v^{2}}
\notag \\
&&-2uvg\frac{\partial ^{2}}{\partial u\partial v}-2ug\frac{\partial }{%
\partial u}-2vg\frac{\partial }{\partial v}.  \label{eq 36}
\end{eqnarray}

On account of (\ref{eq 4}) we get

\begin{equation}
\Delta ^{III}x_{1}=\Delta ^{III}u=\lambda _{11}u+\lambda _{12}v+\lambda
_{13}\left( \frac{a}{2}u^{2}+\frac{b}{2}v^{2}\right) ,  \label{eq 37}
\end{equation}

\begin{equation}
\Delta ^{III}x_{2}=\Delta ^{III}v=\lambda _{21}u+\lambda _{22}v+\lambda
_{23}\left( \frac{a}{2}u^{2}+\frac{b}{2}v^{2}\right) ,  \label{eq 38}
\end{equation}

\begin{equation*}
\Delta ^{III}x_{3}=\Delta ^{III}\sqrt{\omega }=\lambda _{31}u+\lambda
_{32}v+\lambda _{33}\left( \frac{a}{2}u^{2}+\frac{b}{2}v^{2}\right) .
\end{equation*}

Applying (\ref{eq 36}) on the coordinate functions $x_{i},i=1,2$ of the
position vector $\boldsymbol{x}$ and by virtue of (\ref{eq 37}) and (\ref{eq
38}) we find respectively

\begin{equation}
\Delta ^{III}u=-2ug=\lambda _{11}u+\lambda _{12}v+\lambda _{13}\left( \frac{a%
}{2}u^{2}+\frac{b}{2}v^{2}\right) ,  \label{eq 39}
\end{equation}

\begin{equation}
\Delta ^{III}v=-2vg=\lambda _{21}u+\lambda _{22}v+\lambda _{23}\left( \frac{a%
}{2}u^{2}+\frac{b}{2}v^{2}\right) .  \label{eq 40}
\end{equation}

Putting $v=0$ in (\ref{eq 39}), we obtain that

\begin{equation*}
-2a^{2}u^{3}-2u=\lambda _{11}u+\lambda _{13}\frac{a}{2}u^{2}.
\end{equation*}

This implies that $a$ must be zero$.\ $Putting $u=0$ in (\ref{eq 40}), we
obtain that

\begin{equation*}
-2b^{2}v^{3}-2v=\lambda _{22}v+\lambda _{23}\frac{b}{2}v^{2}.
\end{equation*}

This implies that $b$ must be zero$,$ which is clearly impossible, since $%
a>0 $ and$\ b>0.$

\end{document}